\documentclass{article}[12pt]
\usepackage{amsmath,amsfonts,amsthm}
\usepackage{amssymb}

\usepackage{epsfig}

\usepackage{graphics}
\usepackage{graphicx}

\usepackage{latexsym}
\usepackage{graphics}
\usepackage{amsmath}
\usepackage{amsthm}
\usepackage{xspace}
\usepackage{amssymb}
\usepackage{color}
\usepackage{cite}
\usepackage{latexsym}
\usepackage{graphics}
\usepackage{amsmath}
\usepackage{amsthm}
\usepackage{xspace}
\usepackage{amssymb}
\usepackage{epsfig}

\newcommand{\beql}[1]{\begin{equation}\label{#1}}
\newcommand{\eeql}{\end{equation}}
\newcommand{\eqn}[1]{(\ref{#1})}

\newcommand{\pr}{\mathbb{P}}
\newcommand{\E}{\mathbb{E}}

\newcommand{\cx}{{\cal X}}

\newcommand{\bI}{{\bf I}}

\newcommand{\xone}{^{(1)} \! x}
\newcommand{\xtwo}{^{(2)} \! x}

\newcommand{\xind}{^{(\ell)} \! x}

\newcommand{\odin}{^{(1)} \! }
\newcommand{\dva}{^{(2)} \! }
\newcommand{\tri}{^{(3)} \! }
\newcommand{\chet}{^{(4)} \! }
\newcommand{\ind}{^{(\ell)} \! }

\newtheorem{thm}{Theorem}
\newtheorem{lem}[thm]{Lemma}

\newtheorem{cor}[thm]{Corollary}

\newtheorem{definition}[thm]{Definition}

\begin{document}

\title{Large-scale Join-Idle-Queue system with general service times 
}

\author
{
Sergey Foss \\
Heriot-Watt University\\
EH14 4AS Edinburgh, UK\\ and\\
 Novosibirsk State University\\
\texttt{s.foss@hw.ac.uk}
\and
Alexander L. Stolyar \\
University of Illinois at Urbana-Champaign\\
Urbana, IL 61801, USA \\
\texttt{stolyar@illinois.edu}
}

\date{\today}

\maketitle

\begin{abstract}

A parallel server system with $n$ identical servers is considered. The service time distribution has 
a finite mean $1/\mu$, but otherwise is arbitrary. Arriving customers are be routed to one of the servers
immediately upon arrival.
Join-Idle-Queue routing algorithm is studied, under which an arriving customer is sent to an idle server, if such is available, and to a randomly uniformly chosen server, otherwise.
We consider the asymptotic regime where $n\to\infty$ and the customer input flow rate is $\lambda n$.
Under the condition $\lambda/\mu<1/2$, we prove that, as $n\to\infty$, the sequence of (appropriately scaled) 
stationary distributions concentrates at the natural equilibrium point, with the fraction of occupied servers
being constant equal $\lambda/\mu$. In particular, this implies that the steady-state probability of an arriving customer waiting for service vanishes.

\end{abstract}




\section{Introduction}
\label{sec-intro}

We consider a parallel server system consisting of $n$ servers, processing
a single input flow of customers. 
The service time of any customer by any server has the same distribution with finite mean $1/\mu$.
Each customer has to be assigned (routed) to one of
the servers immediately upon arrival. (This model is sometimes referred to as ``supermarket'' model.)
We study a Join-Idle-Queue routing algorithm, under which an arriving customer is sent to an idle server, if such is available; if there are no idle servers, a customer is sent to one of the servers chosen uniformly at random.

We consider an asymptotic regime such that $n\to\infty$ and the input rate is $\lambda n$, where the system load $\lambda/\mu<1$. Thus, the system remains subcritically loaded. Under the additional assumption
that the service time distribution has {\em decreasing hazard rate} (DHR), it is shown in \cite{St2014_pull}
that the following property holds.

{\em Asymptotic optimality: As $n\to\infty$, the sequence of the system stationary distributions is such that the fraction 
of occupied servers converges to constant $\lambda/\mu$; consequently, the steady-state probability of an 
arriving customer being routed to a non-idle server vanishes.}

The results of \cite{St2014_pull} apply to far more general 
systems, where servers may be non-identical. However, the analysis in \cite{St2014_pull} does rely in essential way on the DHR assumption on the service times; under this assumption the system process has {\em monotonicity} property, which is a powerful tool for analysis. Informally speaking, monotonicity means that two versions of the process, such that the initial state of the first one is dominated (in the sense of some natural partial order) by that of the second one,
can be coupled so that this dominance persists at all times.

When the service time distribution is general, the monotonicity under JIQ no longer holds, which requires a different
approach to the analysis. In this paper we prove the following

{\bf Main result} (Theorem~\ref{th-main} in Section~\ref{sec-gen-service}): {\em The asymptotic optimality holds for 
an arbitrary service time distribution, if the system load $\lambda/\mu< 1/2$.}

We believe that condition $\lambda/\mu<1/2$ is purely technical
(required for the proof in this paper)
and that our main result in fact 
holds for $\lambda/\mu<1$, i.e. as long as the system is stable. 
This will be discussed in more detail in Section~\ref{subsec-intuit}.

The key feature of the JIQ algorithm (as well as more general {\em pull-based} algorithms 
\cite{BB08,G11,St2014_pull,St2015_pull}), is that it does not utilize any information about the current state 
of the servers besides them being idle or not. This allows for a very efficient practical implementation,
requiring very small communication overhead between the servers and the router(s) 
\cite{G11,St2014_pull,St2015_pull}. In fact, in the asymptotic regime that we consider, JIQ is much superior 
to the celebrated ``power-of-d-choices'' (or Join-Shortest-Queue(d), or JSQ(d)) algorithm \cite{VDK96,Mitz2001,BLP2012-jsq-asymp-indep,BLP2013-jsq-asymp-tail},
 in terms of both performance and communication overhead (see \cite{St2014_pull,St2015_pull}
for a detailed comparison). The JSQ(d) algorithm routes a customer to the shortest queues among
the $d \ge 1$ servers picked uniformly at random.

We note that when the service time distribution is general, there is
no monotonicity under JSQ(d) (just like under JIQ in our case), and this also makes the analysis far more 
difficult. Specifically, the result for JSQ(d), which is a counterpart of our main result for JIQ,
is Theorem 2.3 in \cite{BLP2012-jsq-asymp-indep}, which shows the asymptotic independence of
individual server states. (Our main result also implies asymptotic independence of server 
states; see formal statement in Corollary~\ref{cor-asym-indep}.)
Theorem 2.3 in \cite{BLP2012-jsq-asymp-indep} imposes even stronger assumptions than ours,
namely a finite second moment of the service time and load $\lambda/\mu < 1/4$ (for non-trivial 
values of $d$, which are $d\ge 2$); our Theorem~\ref{th-main} only requires a finite first moment
of the service time and load $\lambda/\mu < 1/2$.

In a different asymptotic regime, so called Halfin-Whitt regime (when the system capacity exceeds its load 
by $O(\sqrt{n})$, as opposed to $O(n)$), and Markov assumptions (Poisson input flows and
exponentially distributed service times), JIQ has been recently analyzed in \cite{EsGam2015,MBLW2015}.
These papers study diffusion limits of the system transient behavior; Markov assumptions appear to be essential
for the analysis. Finally, we mention a recent paper \cite{St2015_grand-het}, which proposes and studies 
a version of JIQ for systems with {\em packing constraints} at the servers.

{\bf Paper organization.} Section~\ref{sec-gen-service} gives the formal model and main result,
with Section~\ref{subsec-intuit} discussing the role of condition $\lambda/\mu<1/2$.
A uniform stochastic upper bound on the individual server workload in steady-state is derived 
in Section~\ref{sec-uniform-bound}. Properties of the the process fluid limits are established in 
Section~\ref{sec-fluid-lim}. Section~\ref{sec-main-proof} contains the proof of the main result,
which relies on the above upper bound and fluid limit properties.
Generalizations of the main result are presented in Section~\ref{sec-general}.



{\bf Basic notation.} 
The following abbreviations are used to qualify a convergence of functions: {\em u.o.c.} means {\em uniform on compact sets}, {\em p.o.c.} means {\em convergence at points of continuity of the limit}, and
{\em a.e.} means {\em almost everywhere w.r.t. Lebesgue measure}.
We say that a function is RCLL if it is {\em right-continuous with left-limits}.
A scalar function $f(t), ~t\ge 0$, we will call {\em Lipschitz above} if there exist a constant $L>0$ such that
$f(t_2) - f(t_1) \le L(t_2-t_1)$ for any $t_1 \le t_2$. The norm of a scalar function is $\|f(\cdot)\| \doteq \sup_w |f(w)|$. Inequalities applied to vectors [resp.  functions] are understood componentwise 
[resp. for every value of the argument]. Symbol $\Rightarrow$ signifies convergence of random elements in distribution. Indicator of event or condition $B$ is denoted by $\bI(B)$.
Abbreviation WLOG means {\em without loss of generality}.

\section{Model and main result}
\label{sec-gen-service}

We consider a service system, consisting of $n$ parallel servers. The system is homogeneous in that all servers are identical, with the same customer service time distribution, given by the cdf $F(\xi), \xi\ge 0$.
This distribution has finite mean, which WLOG can be assumed to be $1$:
$$
\int_0^\infty F^c(\xi) = 1, ~~~\mbox{where}~F^c(\xi) \doteq 1 - F(\xi).
$$
Otherwise, the cdf $F(\cdot)$ is arbitrary. The service/queueing discipline at each server is arbitrary, as long as it is work-conserving and non-idling.

Customers arrive as a Poisson process. (This assumption can be relaxed to a renewal arrival process;
see Section~\ref{sec-general}.) The arrival rate is $\lambda n$, where $\lambda<1$, so that the system load is strictly subcritical.

The routing algorithm is Join-Idle-Queue (JIQ), which is defined as follows. (The JIQ algorithm can be viewed, in particular, as a specialization of the PULL algorithm \cite{St2014_pull,St2015_pull} to a homogeneous system with ``single router.'')

\begin{definition}[JIQ]
\label{def-pull-basic1}
An arriving customer is routed to an idle server, if there is one available. Otherwise, it is routed to server chosen uniformly at random.
\end{definition}

We consider the sequence of systems with $n\to\infty$.
From now on, the upper index $n$ of a variable/quantity will indicate that it pertains to the system with $n$ servers,
or $n$-th system.
Let $W_i^n(t)$ denote the workload, i.e. unfinished work, in queue $i$ at time $t$ in the $n$-th system.
Consider the following {\em fluid-scaled} quantities:
\beql{eq-x-def}
x^n_w(t) \doteq (1/n) \sum_i \bI\{W_i^n(t)> w\}, ~~ w \ge 0.
\end{equation}
That is, $x^n_w(t)$ is the fraction of servers $i$ with $W_i^n(t)> w$.
Then
 $x^n(t)=(x^n_w(t), ~w\ge 0)$
is the system state at time $t$;
$\rho^n(t) \equiv  x^n_0(t)$ is the fraction of busy servers (the instantaneous system load).

For any $n$, the state space of the process $(x^n(t), ~t\ge 0)$ is a subset of a common (for all $n$) 
state space $\cx$, whose elements
$x=(x_w, ~w\ge 0)$ are non-increasing RCLL functions of $w$, 
with values $x_w\in [0,1]$. 
This state space $\cx$ is equipped with Skorohod metric, topology and corresponding Borel $\sigma$-algebra. 

Then, for any $n$, process $x^n(t), ~t\ge 0$ is Markov with state space $\cx$,
and sample paths being RCLL functions (with values in $\cx$), which are in turn elements of (another) Skorohod
space. (The Skorohod spaces that we defined play no essential role in our analysis; we need
to specify them merely to make the process well-defined.)

Stability (positive Harris recurrence) of the process $(x^n(t), ~t\ge 0)$, for any $n$, is straightforward to verify.
Indeed, as long as a server remains busy, it receives each new arrival with probability at most $1/n$, 
and therefore receives the new work at the average rate at most $\lambda$. (We omit the details of stability proof.)
Thus, the process has unique stationary distribution. Let $x^n(\infty)$ be a random element whose distribution is the stationary distribution of the process; in other words,
this is a random system state in stationary regime.

The system {\em equilibrium point} $x^*  \in \cx$ is defined as follows.
Let $\Phi^c(w)$ denote the complementary (or, tail) distribution function of the steady-state residual service time;
the latter is the steady-state 
residual time of a renewal process with renewal time distribution function $F(\cdot)$. We have
$$
\Phi^c(w) = \int_w^\infty F^c(\xi) d\xi, ~~w\ge 0.
$$
Then, 
$$
x^* = (x^*_w = \lambda \Phi^c(w), ~w\ge 0) \in \cx.
$$
In particular, the equilibrium point is such that ``the fraction of occupied servers'' $x^*_0 = \lambda$.
Our main result is the following

\begin{thm}
\label{th-main}
If $\lambda<1/2$, then $x^n(\infty) \Rightarrow x^*$ as $n\to\infty$. 
\end{thm}

The theorem shows, in particular, that  if $\lambda < 1/2$, then as $n\to\infty$ the steady-state 
probability of an arriving customer waiting for service (or sharing a server with other customers) vanishes.
Theorem~\ref{th-main} easily generalizes to the case when: (a) arrival process is renewal, (b) some or all servers
may have finite buffers, and (c) there may be some bias in the routing when all servers are busy.
(These generalizations are described in Section~\ref{sec-general}.)

Theorem~\ref{th-main} implies the following 

\begin{cor}
\label{cor-asym-indep}
Assume $\lambda<1/2$. 
Suppose
that
JIQ is completely symmetric with respect to the servers. Specifically, if at the time of a customer arrival there are idle servers, the customer is routed to one of them chosen uniformly at random. Then
the states of individual servers in stationary regime are 
asymptotically independent. Moreover, 
for any fixed $m$, the stationary distribution of $(W_1^n, \ldots, W_m^n)$ converges to
that of $(\widetilde W_1, \ldots, \widetilde W_m)$, with i.i.d. components such that 
$\pr\{\widetilde W_1 > w\} = x^*_w = \lambda \Phi^c(w), ~w\ge 0$.
\end{cor}

Indeed, by symmetry with respect to servers, the stationary distribution of $(W_1^n, \ldots, W_m^n)$,
i.e. of the residual work on the fixed set of servers $1, \ldots, m$, is same as that on 
a set of $m$ servers, {\em chosen uniformly at random.} But, $x^n(\infty)$, which describes the overall distribution of server
workloads in the system, converges in distribution to the non-random point $x^*$. This implies Corollary~\ref{cor-asym-indep}.

\subsection{Discussion of condition $\lambda < 1/2$.}
\label{subsec-intuit}

The approach we use to establish the convergence of stationary distributions in Theorem~\ref{th-main} is 
as follows. We find a set $A \in \cx$ and a fixed finite time $T$, such that, with high probability, for all large $n$, (a) $x^n(\infty)\in A$ and (b) $x^n(0)\in A$ implies that $x^n(T)$ is close to $x^*$. 
Property (b) is key.
When $n$ is large, the trajectory  $x^n(t)$ is ``almost deterministic.'' (In fact, 
the problem reduces to the analysis of ``fluid limit'' trajectories, which are the limits of 
$x^n(t)$ as $n\to\infty$.) Then, informally speaking, property (b) above reduces to the property (b'):
trajectories $x^n(t)$ converge to $x^*$ as $t\to\infty$.
The absence of process monotonicity (described in Section~\ref{sec-intro}) makes proving (b') difficult. 
We now describe -- very informally -- the key idea, which we use in our proof of convergence (b'), and which relies on the condition $\lambda<1/2$.

Suppose $n$ is large. 
Consider an initial state $x^n(0)$, such that {\em the total amount of (fluid-scaled, i.e. multiplied by $1/n$) unfinished work is upper bounded by $C<\infty$}.
Pick $\alpha$ such that $\alpha > \lambda$
and $\alpha+\lambda<1$; this can be done if and only if $\lambda<1/2$.
Then, at some finite time $\tau$, the system must reach a state with $\alpha n$ servers being idle. (Otherwise, if at least $(1-\alpha)n$ servers would continue to be busy as time goes to infinity, the unfinished work would become negative, since $1-\alpha > \lambda$.) Denote by $S_\alpha$ the set of those $\alpha n$ servers, which are idle at time $\tau$.
Starting time $\tau$, WLOG, assume that all new arriving customers go to an idle server in $S_\alpha$,
as long as there is one available. Consider the subsystem, consisting only of the servers in $S_\alpha$;
starting time $\tau$ and until the (random) time when {\em all} servers in $S_\alpha$ become busy,
the behavior of this subsystem is obviously equivalent to that of the infinite-server system, $M/GI/\infty$,
with idle initial state. If $n$ is large, the behavior of $x^n(t)$ {\em for such $M/GI/\infty$ system} is ``almost deterministic''
and such that the (scaled) number of occupied servers $x_0^n(t)$ in it is ``almost monotone increasing, 
converging to $\lambda  < \alpha$" and, moreover, $x^n(t)$ ``converges" to $x^*$. But this means that after time $\tau$ the subsystem $S_\alpha$ will ``always'' have idle servers, which in turn means that {\em its} state will ``converge'' to $x^*$ as $t\to\infty$. 
Also, after time $\tau$,
the subsystem consisting of the servers outside $S_\alpha$ will ``never'' receive any new arrivals and will ``eventually'' empty. Thus, $x^n(t)$ for our entire system ``converges'' to $x^*$.

Turning the key intuition, described above informally, into a formal proof is the subject of the rest of this paper.
Set $A \in \cx$ is picked by using a constructed uniform in $n$ upper bound on the stationary distribution of the workload of an individual server. 
The states in $A$ are such that the total (scaled) workload is not necessarily upper bounded by a constant $C$
(in fact, if the second moment of the service time is infinite,
 the steady-state total workload in the system is infinite with probability $1$);
however, for states in $A$ the (scaled) workload is bounded by $C$ on a close-to-$1$ fraction of servers -- this suffices for the proofs.
The property (b') is proved uniformly for fluid limits starting from $A$ -- from here we obtain that (b) holds for the pre-limit processes with high probability, uniformly for all large $n$.

As explained above, our proof of Theorem~\ref{th-main} relies in essential way on condition 
$\lambda<1/2$. However, we believe that this condition is purely technical, and Theorem~\ref{th-main} in fact holds for any $\lambda < 1$. Establishing this fact will most likely require a different proof approach, although some elements of the analysis in this paper may turn out to be useful for the proof of a more general result.

\section{Uniform upper bound on a server workload distribution}
\label{sec-uniform-bound}

Throughout this section, consider a fixed $\lambda<1$. Consider an M/GI/1 system, with arrival rate $\lambda$
and service time distribution $F(\cdot)$. Let us view its workload process as regenerative with renewal points being time instants when a customer arrives into idle system. For each $w\ge 0$,
denote by $x_w^{**}$ 
the expectation of the total time during one renewal cycle when the workload is greater than $w$. Clearly, 
$x_w^{**}$ is non-increasing, $x_0^{**} = 1/(1-\lambda)$ (the expected busy period duration) and 
$x_w^{**} \to 0, ~w\to\infty$. (We will not use the exact value of $x_0^{**}$. Also,
$x_w^{**}$ is continuous in $w$, but we will not use this fact either.)

Now consider our system with any fixed $n$.
Consider a specific server $i$. Consider our Markov process sampled at the ``renewal'' instants when there
is an arrival into idle server $i$. 
Time intervals between the ``renewal'' instants are ``renewal cycles''. Of course, such ``renewal cycles'' are not i.i.d., the law of a cycle depends on the state of the entire system at the renewal point from which the cycle starts. However, there are uniform bounds that apply to any cycle. For a fixed $w\ge 0$,
the expected total time within one cycle when $W^n_i>w$, is upper bounded by $x_w^{**}$;
indeed, as long as the server remains busy, the probability that a new arrival will be routed to it is at most $1/n$
(either $1/n$ or $0$); therefore, as long as the server remains busy, the instantaneous arrival rate into it, is upper bounded by $(\lambda n) / n = \lambda$. The mean duration of each cycle is lower bounded by the 
mean service time of one customer, i.e. by $1$. Therefore,
\beql{eq-m-g-1-bound}
\pr\{W^n_i(\infty) > w\} \le x_w^{**}, ~~ w\ge 0,
\end{equation}
where, recall, $x_w^{**} \to 0, ~w\to\infty$.
Bound \eqn{eq-m-g-1-bound} implies
the following fact.

\begin{lem}
\label{lem-m-g-1-bound}
Let $\lambda<1$. Then, for any $n$,
$
\E x^n_w(\infty) \le x^{**}_w, ~~w\ge 0.
$
\end{lem}

\section{Fluid limits}
\label{sec-fluid-lim}

In this section we introduce different types of the process fluid limits, which will be used later in the analysis.

We will assume that, given a fixed initial state $x^n(0)$, the realization of the process is determined by a common (for all $n$) set of driving processes. Specifically, there is a common, rate $1$, Poisson process, $\Pi(t), ~t\ge 0,$;
the number of arrivals in the $n$-th system by time $t$ is $\Pi(n\lambda t)$. There is also a common sequence
of i.i.d. random variables with distribution $F(\cdot)$, which determines the service times of arriving customers
(in the order of arrivals). 
Let $G^n(t,w), ~t\ge 0, ~w \ge 0,$ be the number of customer arrivals 
in the $n$-th system, by time $t$, with the service times greater than $w$.
Let $g^n(t,w) = (1/n) G^n(t,w)$ and $g(t,w) \doteq \lambda t F^c(w)$.
 Then, we have the following functional strong law of large numbers (FSLLN):
\beql{eq-lln-driving}
\| g^n(t,\cdot) - g(t,\cdot)\| \to 0, ~\mbox{as}~n\to\infty,~~~u.o.c.~\mbox{(in $t$)}, ~~~w.p.1.
\end{equation}
Indeed, for any fixed $t>0$, the total number of arrivals in $[0,t]$, scaled by $1/n$, converges to $\lambda t$ w.p.1;
this and Glivenko-Cantelli theorem (cf. \cite{Billingsley-95}, Theorem 20.6, page 269) imply that 
$\| g^n(t,\cdot) - g(t,\cdot)\| \to 0$, w.p.1. But, all $g^n(t,w)$ and $g^n(t,w)$ are non-decreasing in $t$, and $g(t,w)$ 
is continuous in $t$; this easily implies that the convergence in \eqn{eq-lln-driving} is uniform w.p.1.

The routing of arriving customers to idle servers, when such are available, is completely arbitrary WLOG;
it will be specified later, in a way convenient for the analysis. The routing of arriving customers to the servers, in cases when all servers are busy is determined by a sequence if i.i.d. random variables, uniformly distributed in
$[0,1)$; these random variables are used sequentially ``as needed''; in the $n$-th system, a customer is routed
to server $i$ if the corresponding random variable value is in $[(i-1)/n,i/n)$. (The specific construction of routing to busy servers will not be important; we need to specify it somehow, to have the process well defined.)

It will be convenient for every $n$, in addition to the actual system with $n$ servers, to consider the corresponding infinite server
system; in such system all arrivals always go to idle servers.
For a given $n$, for the infinite server system the fluid-scaled quantities $x^n_w(t)$
are still defined by \eqn{eq-x-def}, i.e.
as the total number of servers with $W_i^n > w$, multiplied by $1/n$. 

For every $t\ge 0$, let us define $x^\uparrow(t) = (x_w^\uparrow(t), w\ge 0)\in \cx$,
$$
x_w^\uparrow(t) = \int_0^t F^c(w+t-\theta) \lambda d\theta = \int_w^{w+t} F^c(\xi) \lambda d\xi.
$$
Clearly, $x_w^\uparrow(t)$ is non-decreasing in $t$, and 
$$
x_w^\uparrow(t) \uparrow x^*_w, ~~t\to\infty, ~~~\forall w.
$$
As functions of $w$, all $x_w^\uparrow(t)$ and $x^*_w$ are non-negative, continuous, non-decreasing
and converging to $0$ as $w\to \infty$; therefore, the above pointwise convergence 
implies uniform convergence
$$
\|x^\uparrow(t) - x^* \| \to 0, ~~t\to\infty.
$$
The following Lemma~\ref{lem-infinite} is a standard fact.
Informally speaking, it states that $x^\uparrow(\cdot)$ is the ``fluid limit'', in $n\to\infty$, 
of $x^n(\cdot)$ for the infinite-server system, starting from idle initial state.
We state this fact in a form that is convenient for our analysis, and since it easily follows from
\eqn{eq-lln-driving}, we give a proof as well.

\begin{lem}
\label{lem-infinite}
Fix arbitrary $\lambda\ge 0$. (Here  $\lambda \ge 1$ is allowed.) 
Let $x^n(\cdot)$ be the process describing the infinite-server system, starting from idle initial state,
that is, $x^n_0(0)=0$. Then, w.p.1, 
\beql{eq-infinite-fl}
\| x^n(t) - x^\uparrow(t)\| \to 0, ~~~u.o.c.
\end{equation}
\end{lem}

{\em Proof.}  Fix $t$ and $w$. By definition, $x^n_w(t)$ is the scaled number of customers 
in the system, having the residual service time greater than $w$. 
A customer arriving at time $\theta\in [0,t]$ counts into that number if and only if its service 
time is greater than $t+w-\theta$. Let points $0=t_0 < t_1 < \ldots <t_\kappa=t$ partition
the interval $[0,t)$ into $\kappa$ subintervals $[t_k,t_{k+1})$. (W.p.1 there are no arrivals at $t$.)
Then,
$$
\sum_{k=0}^{\kappa-1} [g^n(t_{k+1}, t +w - t_{k}) - g^n(t_{k}, t +w - t_{k})] \le 
x^n_w(t) \le \sum_{k=0}^{\kappa-1} [g^n(t_{k+1}, t +w - t_{k+1}) - g^n(t_{k}, t +w - t_{k+1})].
$$
By \eqn{eq-lln-driving}, w.p.1 the lower and upper bounds converge to 
$$
\sum_{k=0}^{\kappa-1} \lambda [t_{k+1}-t_k] F^c(t +w - t_{k}) ~~\mbox{and}~~
\sum_{k=0}^{\kappa-1} \lambda [t_{k+1}-t_k] F^c(t +w - t_{k+1}),
$$
respectively. Considering a sequence of partitions with maximum subinterval size vanishing,
and taking into account that $F^c$ is non-increasing, we obtain probability 1 convergence
$x^n_w(t) \to x_w^\uparrow(t)$. Since $x_w^\uparrow(t)$ and all $x^n_w(t)$ are non-negative non-increasing 
in $w$, $x_w^\uparrow(t)$ is continuous in $w$,
and $x_w^\uparrow(t) \to 0$ as $w\to\infty$, we obtain probability 1 convergence
$\|x^n(t)  - x^\uparrow(t)\| \to 0$, for any $t$; since $x_w^\uparrow(t)$ is continuous non-decreasing in $t$,
this convergence is u.o.c. in $t$.
$\Box$

Sometimes, it will be convenient to divide the set of servers into two or more subsets, and keep track 
of the workloads in those subsets separately. 
For example, suppose at time $0$ the set of all servers, let us call it $S$,
is divided (for each $n$) at time $0$ into two non-intersecting subsets, $S_1$  and 
$S_2$, and these subsets do not change with time.
Then, for $\ell=1,2$, $\xind^n_w(t)$ is the fraction of servers (out of the total number $n$)
which are in $S_\ell$ and have workload $W^n_i> w$, $w \ge 0$; $\ind\rho^n(t)= ~ \xind^n_0(t)$.
Of course, $x(t) = ~ \xone(t) + ~\xtwo(t)$. However, often we will consider $\xind(t)$ for only one 
of the subsets $S_\ell$.

The following fact is a corollary of Lemma~\ref{lem-infinite}.

\begin{lem}
\label{lem-finite}
Let $0\le \lambda < 1$ and let $\lambda<\alpha <1$. Consider the finite server system. Assume that 
for all $n$, the initial states are such that $x^n_0(0)=\rho^n(0) \le 1- \alpha$.
For each $n$, consider the subset $S_1 = S_1(n)$, consisting of $\alpha n$ servers that are 
initially idle. 
Assume WLOG that any new arrival will go to an idle server in $S_1$, if there is one available.
Then, w.p.1, 
the following holds:
\beql{eq-finite-fl}
\| \odin x^n(t) - x^\uparrow(t) \| \to 0, ~~~u.o.c.,
\end{equation}
and for any fixed $t$, for all sufficiently large $n$, all new arrivals in $[0,t]$ will go to idle servers
in $S_1$.
\end{lem}

{\em Proof.}  
The behavior of the system restricted to subset $S_1$ of servers is equivalent to that of the infinite server system starting from idle state, {\em as long as there are idle servers in $S_1$}. 
By Lemma~\ref{lem-infinite}, w.p.1 the trajectory of the (scaled) infinite-server system converges (u.o.c.)
to the trajectory $x^{\uparrow}(t)$, such that the (scaled) number of occupied server increases and converges to $\lambda<\alpha$. This implies that w.p.1. the following holds for the system restricted to subset $S_1$:  
for any fixed time $t\ge 0$, for all sufficiently large $n$, subset $S_1$ will have idle servers in the 
entire interval $[0,t]$, and then the system behavior coincides with that of the infinite-server system.
This property implies \eqn{eq-finite-fl}, and contains the last statement of the lemma.
$\Box$

Let $\ind W^n(t)$ denote the total (fluid-scaled) unfinished
work at time $t$ within a given subset $S_\ell$ of servers:
$$
\ind W^n(t) = \int_0^\infty ~\ind x^n_w(t) dw.
$$
The case $S_\ell = S$ is allowed.

Denote by $\ind W^{a,n}(t)$ and $\ind W^{d,n}(t)$ the amount of (fluid-scaled) 
work that, respectively, arrived into and processed by subset $S_\ell$ in the interval $[0,t]$.
Denote by $\ind \rho^{a,n}(t)$ the (fluid-scaled) number of arrivals in $[0,t]$ into $S_\ell$,
that went into idle servers; such arrivals, and only they, cause $+1/n$ jumps of $\ind \rho^{n}$.
Analogously, let $\ind \rho^{d,n}(t)$ denote the (fluid-scaled) number of times in $[0,t]$ 
when a customer service completion occurred in $S_\ell$,
that left a server idle; such departures, and only they, cause $-1/n$ jumps of $\ind \rho^{n}$.
Functions $\ind W^{a,n}(t)$, $\ind W^{d,n}(t)$, $\ind \rho^{a,n}(t)$ and $\ind \rho^{d,n}(t)$
are non-decreasing by definition, equal to $0$ at $t=0$.
The following relations obviously hold for all $t\ge 0$:
\beql{eq-conserv1}
\ind W^n(t) = ~\ind W^{a,n}(t) - ~ \ind W^{d,n}(t), ~~
\ind \rho^n(t) = ~\ind \rho^{a,n}(t) - ~ \ind \rho^{d,n}(t),
\end{equation}
\beql{eq-conserv2}
\ind W^{d,n}(t) = \int_0^t ~\ind \rho^n(\xi) d\xi.
\end{equation}
For future reference let us also note the obvious fact that {\em if there were no new arrivals into
$S_\ell$ in some time interval $(t_1,t_2]$, then}
\beql{eq-conserv3}
\ind W^n(t_2) - ~ \ind W^n(t_1) = - (~\ind W^{d,n}(t_2) - ~\ind W^{d,n}(t_1)) = 
- \int_{t_1}^{t_2} ~\ind \rho^n(\xi) d\xi.
\end{equation}

\begin{lem}
\label{lem-workloads}
Let $\lambda \ge 0$. Consider the finite server system. 
For each $n$ consider a subset $S_1=S_1(n)$, consisting of $\sigma n$ servers, $0\le \sigma \le 1$.
(The case $\sigma=1$ is when $S_1=S$.)
Consider a fixed sequence (in $n$) of initial states, such that $\odin W^n(0) \le C < \infty, ~\forall n$.
Then, w.p.1, for any subsequence of $n$, there exists a further subsequence, along which
the following holds:
\beql{eq-W-conv}
\odin W^n(t) \to ~\odin W(t), ~~~u.o.c.,
\end{equation}
where $\odin W(\cdot)$ is a Lipschitz continuous function with $\odin W(0) \le C$;
\beql{eq-rho-conv}
\odin \rho^n(t) \to ~\odin \rho(t), ~~~p.o.c.,
\end{equation}
where $\odin \rho(\cdot)$ is a RCLL function, which is Lipschitz above and 
$\odin \rho(t) \in [0,\sigma], ~\forall t$;
\beql{eq-W-deriv111}
\odin W'(t) \le \lambda - ~\odin\rho(t), ~~a.e.
\end{equation}
\end{lem}

{\em Proof.} Within this proof, when we say that a function is Lipschitz continuous (resp., Lipschitz above),
we always mean that it is Lipschitz continuous (resp., Lipschitz above) {\em uniformly in $n$.}

From FSLLN \eqn{eq-lln-driving} we have the following fact. W.p.1, for any subsequence of $n$, there exists a further subsequence, along which
$$
\odin \rho^{a,n}(t) \to ~\odin \rho^{a}(t), ~~\odin W^{a,n}(t) \to ~\odin W^{a}(t), ~~u.o.c.,~~~n\to\infty.
$$
where $\odin \rho^{a}(\cdot)$ and $\odin W^{a}(\cdot)$ are Lipschitz continuous non-decreasing,
with Lipschitz constant equal $\lambda$.
Also, clearly, all functions $\odin W^{d,n}(\cdot)$ are non-decreasing Lipschitz continuous,
so that we can choose a further subsequence, if necessary, along which
$$
\odin W^{d,n}(t) \to ~\odin W^{d}(t), ~~u.o.c.,~~~n\to\infty,
$$
where $\odin W^{d}(\cdot)$ is Lipschitz continuous non-decreasing.
This implies \eqn{eq-W-conv} with $\odin W(t) = ~\odin W^{a}(t) - ~\odin W^{d}(t)$.

To show \eqn{eq-rho-conv}, observe that non-decreasing functions $\odin \rho^{d,n}(t)$ are uniformly
bounded on any finite interval (because functions $\odin \rho^{a,n}(t)$ and $\odin \rho^{n}(t)$
are, along the chosen subsequence).
Then, we can choose a further subsequence, if necessary, such that
\beql{eq-777}
\odin \rho^{d,n}(t) \to ~\odin \rho^{d}(t), ~~p.o.c.,~~~n\to\infty,
\end{equation}
where $\odin \rho^{d}(\cdot)$ is RCLL non-decreasing. 
(Here we use a version of Helly's selection theorem; cf. \cite{Billingsley-95}, Theorem 25.9, page 336.)
This proves 
\eqn{eq-rho-conv} with $\odin \rho(t) = ~\odin \rho^{a}(t) - ~\odin \rho^{d}(t)$.

Note that the {\em p.o.c.} convergence in \eqn{eq-777} implies {\em a.e.} (in $t$) convergence.
Then, by taking limit in \eqn{eq-conserv2}, we obtain
$$
\odin W^{d}(t) = \int_0^t ~\odin \rho(\xi) d\xi.
$$
This and the fact that  $\odin W^{a}(\cdot)$ is Lipschitz continuous 
with Lipschitz constant $\lambda$, imply \eqn{eq-W-deriv111}.
$\Box$

\section{Proof of Theorem~\ref{th-main}}
\label{sec-main-proof}

Here we only consider the finite systems (with $n$ servers in $n$-th system).
Consider a fixed $\lambda < 1/2$. 

By Lemma~\ref{lem-m-g-1-bound},
for each $n$ we have
$
\E x_w^n(\infty) \le x_w^{**},
$
where $x_w^{**}$ is non-increasing and $\lim_{w\to\infty} x_w^{**} = 0$.
Then for any 
$\delta_1>0$
we can choose a sufficiently large $b$, such that $\E x_b^n(\infty) \le \delta_1$.
This in turn implies that for any $\epsilon>0$ and any $\delta>0$ we can pick sufficiently large $b>0$,
such that
\beql{eq-prob-bound}
\pr\{x_b^n(\infty) \le \delta\} \ge 1-\epsilon, ~~\forall n.
\end{equation}

For each $n$ consider the {\em stationary version of process} $x^n(\cdot)$;
then, for any $t$,  
$x^n(t)$ is equal in distribution to $x^n(\infty)$ (by the definition of the latter). 
Choose $\delta > 0$ small enough so that
$\lambda+\delta < 1/2$.
For this $\delta$ and arbitrarily small fixed $\epsilon>0$, choose $b>0$ such that \eqn{eq-prob-bound} holds.
Then, \eqn{eq-prob-bound} implies
\beql{eq-prob-bound111}
\pr\{\mbox{condition \eqn{eq-cond-n} holds} \} \ge 1-\epsilon, ~~\forall n,
\end{equation}
\beql{eq-cond-n}
\mbox{$x^n(0)$ is such that 
$\exists$ a subset $S_2=S_2(n)$ 
of $(1-\delta)n$ servers, each with workload at most $b$.}
\end{equation}

Then, to complete the proof of Theorem~\ref{th-main}, it suffices to prove the following

\begin{lem}
\label{lem-conv}
For any $\delta_2>0$ there exists $T>0$, which depends on $\epsilon, \delta, b$,
such that, uniformly on fixed initial states $x^n(0)$ satisfying \eqn{eq-cond-n}, 
\beql{eq-conv}
\pr\{ \| x^n(T) - x^*\|   \le \delta + \delta_2   ~|~ x^n(0) \} \to 1, ~~n\to\infty.
\end{equation}
\end{lem}

Indeed, if Lemma~\ref{lem-conv} holds, then for $\delta, \epsilon, b, \delta_2, T$ chosen as 
specified above, and arbitrarily small $\epsilon_2>0$, for all sufficiently large $n$,
uniformly on $x^n(0)$ satisfying \eqn{eq-cond-n},
$$
\pr\{ \| x^n(T) - x^*\|   \le \delta + \delta_2   ~|~ x^n(0) \} \ge 1-\epsilon_2.
$$
This and \eqn{eq-prob-bound111} imply that for all sufficiently large $n$
$$
\pr\{ \| x^n(T) - x^*\|   \le \delta + \delta_2\} \ge (1-\epsilon)(1-\epsilon_2).
$$
But, $\delta, \delta_2, \epsilon, \epsilon_2$ can be chosen arbitrarily small, and recall that
$x^n(T)$ is equal in distribution to $x^n(\infty)$. This proves Theorem~\ref{th-main}.

{\em Proof of Lemma~\ref{lem-conv}.} To establish \eqn{eq-conv} it will suffice to show that
for some fixed $T$ the following holds for any fixed sequence of initial states $x^n(0)$,
satisfying \eqn{eq-cond-n}: the process can be 
constructed in such a way that w.p.1 for all sufficiently large $n$,
\beql{eq-conv111}
\| x^n(T) - x^*\|   \le \delta + \delta_2.
\end{equation}

Fix $\tau > 2 b/(\lambda+\delta/2)$. Fix $T > \tau$. (The choice of $T$ will be specified later.)
For each $n$, at initial time $0$, fix a subset of servers $S_2$ as in
condition \eqn{eq-cond-n}; let $S_1 = S \setminus S_2$ be the complementary subset
of servers -- its size is $\delta n$. Clearly, for each $n$,
$$
\dva W^n(0) \le b, ~~\odin \rho^n(t) \le \delta, ~\forall t.
$$
Consider Markov (stopping) time $\tau^n$, defined as the smallest time $t$ in $[0,\tau]$, such that
$\dva \rho^n(t) \le \lambda+\delta/2$; if there is no such $t$, then $\tau^n=\infty$ by convention.
The construction of the process in $[0,T]$ will be as follows: in the interval $[0,\tau^n]$ it is driven by one 
set of driving processes, and in $(\tau^n,T]$ it is driven by a different, independent set of driving processes
with the same law. (However, these two sets of driving processes are common for all $n$.)
In other words, at time $\tau^n$ the process is ``restarted,'' with the state at
$\tau^n$ serving as initial state and with a new independent set of driving processes. By convention, if 
$\tau^n=\infty$, the process is {\em not} restarted.

We see that w.p.1 for all sufficiently large $n$, 
\beql{eq-tau111}
\tau^n < \tau.
\end{equation}
 Indeed, if we apply 
Lemma~\ref{lem-workloads} to $\dva x^n(t)$ starting time $0$, we see that
any fluid limit $(\dva W(\cdot),~ \dva \rho(\cdot))$ that can arise
is such that $\dva W(0) \le b$ and there exists $t' \le \tau/2$ such that $\dva \rho(t') \le \lambda+\delta/2$.
(Otherwise $\dva W(t)$ would become negative.) This implies \eqn{eq-tau111}.

Similarly we see that w.p.1 for all sufficiently large $n$,
\beql{eq-tau222}
\dva W^n(\tau^n) \le b_1 \doteq b + 2 \lambda \tau.
\end{equation}

Now, consider any fixed sequence of $\tau^n < \tau$ and fixed states at $\tau^n$, satisfying
\eqn{eq-tau111} and \eqn{eq-tau222}. (Recall that starting $\tau^n$, the process is controlled
by a new independent set of driving processes.) Starting time $\tau^n$ we keep the subset $S_1$
as it was, but split $S_2$ into two subsets $S_3$ and $S_4$ as follows: $S_4$
will consists of $(1/2)n$ idle (at $\tau^n$) servers (which exist by \eqn{eq-tau111}),
and $S_3 = S_2 \setminus S_4$ will include the remaining 
$[(1-\delta) - 1/2]n = (1/2 - \delta)n$
servers from $S_2$. Clearly, $\tri W^n(\tau^n) = ~\dva W^n(\tau^n) \le b_1$.
To summarize, starting $\tau^n$, the set of servers $S$ is divided into three
subsets, $S_1$, $S_3$ and $S_4$, with sizes $\delta n$, $(1/2-\delta)n$ and 
$(1/2)n$, respectively. Also, WLOG we assume that starting $\tau^n$ all new arrivals go to
subset $S_4$, as long as there are idle servers in it.
Applying Lemma~\ref{lem-finite}, we obtain that
w.p.1 for all sufficiently large $n$, in the interval $[\tau^n,T]$, all new arrivals go to subset $S_4$.

We now specify the choice of $T$. It has to satisfy two conditions. First, it has to be large enough,
so that for any $t\ge T-\tau$, $\|x^\uparrow(t) - x^*\| \le \delta_2/3$.
Second, it has to be large enough so that 
$$
T-\tau > b_1/(\delta_2/3).
$$
Then applying Lemma~\ref{lem-finite} and \eqn{eq-conserv3}, we obtain that
w.p.1 for all sufficiently large $n$,
$$
\|\chet x^n(T) - x^*\| < \delta_2/2,
$$
$$
\tri \rho^n(T) < \delta_2/2;
$$
this in turn implies \eqn{eq-conv111}.
$\Box$

\section{Generalizations}
\label{sec-general}

\subsection{Renewal arrival process}

The assumption that the arrival process is Poisson is made to simplify the exposition. Our main result,
Theorem~\ref{th-main},
 and the analysis easily generalize to the case when the arrival process is renewal; in the $n$-th system the interarrival times are i.i.d.,  equal in distribution to $A/n$, where $A$ is a positive random variable, $\E A = 1/\lambda$. 
 (Mild assumptions on the interarrival time distribution are needed to make sure that
the process is positive Harris recurrent. For example, it suffices that this distribution has an absolutely
continuous component.)
The common process state space contains an additional scalar variable $u$, which is the residual interarrival time; clearly $u^n(\infty) \Rightarrow 0$ as $n\to\infty$. The more general form of 
Theorem~\ref{th-main} is as follows:

{\em If $\lambda<1/2$, then $(u^n,x^n)(\infty) \Rightarrow (0,x^*)$.}

The construction of the uniform stochastic upper bound on a single server workload generalizes as follows. For each $n$ the arrival process into a server, when it is busy,
is dominated by a renewal process which is the thinned with probability $1/n$ arrival process into the system.
(In other words, as before, the dominating arrival process into a server, as long as the server remains busy,
is such that {\em every} new arrival into the system goes to this server with probability $1/n$.)
The interarrival times of this renewal process are i.i.d. with the distribution equal to that of a random variable $A_n$;
its mean is $\E A_n = 1/\lambda$ for any $n$, but the distribution depends on $n$. 
However, as $n\to\infty$, the distribution of $A_n$ converges to the exponential distribution.
(This is a well known property that a thinned with probability $1/n$ and sped up in time by factor $n$ renewal process
converges to Poisson process. And it is easy to check directly, 
since $A_n$ is a sum of the geometrically distributed, with mean $n$, number of independent instances of $A/n$.)
Then, for arbitrarily small $\delta>0$, there exists
a non-negative random variable $A^\delta$, such that $1/\lambda-\delta \le \E A^\delta < 1/\lambda$,
and the distribution of $A^\delta$ is dominated by that of $A_n$ for all sufficiently large $n$.
(For example, if $\tilde A$ has exponential distribution with mean $1/\lambda$, we can choose 
$A^\delta = ((\tilde A \wedge C)- \epsilon) \vee 0$, where $C>0$ is large, $\epsilon>0$ is small, and $\wedge$, $\vee$
denote minimum and maximum, respectively.)
We fix $\delta>0$ such that $1/\lambda-\delta > 1$, and then $\E A^\delta>1$.
For all large $n$, the renewal arrival process with interarrival times distributed as $A^\delta$ (and then 
the arrival rate $1/\E A^\delta < 1$), dominates (pathwise, using natural coupling)
the arrival process into an individual server, as long as the server
remains busy. 
Therefore, the workload during the busy period under interarrival times $A^\delta$,
dominates 
that under interarrival times $A_n$. 
The rest of the construction of the uniform stochastic upper bound on a single server workload is same.
And after this bound is established, the rest of the proof of the main result remains essentially same
as well, with slight adjustments.

\subsection{Biased routing when all servers busy}

Examination of the proof of Theorem~\ref{th-main} shows that the specific rule -- uniform at random -- for 
routing arriving customers when all servers are busy, is only used to obtain the process stability (positive Harris recurrence) and the uniform stochastic upper bound on a single server workload. In the rest of the proof, this specific rule
is not used; we only use the fact that customers must go to idle servers if there are any. But, for the stability and workload upper bound, it suffices that the arrival rate into a server when it is busy is upper bounded by some 
$\bar\lambda < 1$, not necessarily by $\lambda < 1/2$. This shows that Theorem~\ref{th-main} holds as is,
even if routing when all servers are busy is biased in arbitrary way, as long as the probability that a server
receives an arrival does not exceed $(1/n)(\bar\lambda/\lambda)$ for some $\bar\lambda<1$.

\subsection{Finite buffers}

The main result, Theorem~\ref{th-main}, holds as is if we allow some or all servers to have finite buffers (of same or different sizes). If a server has finite buffer of size $B\ge 1$, and already has $B$ customers, any new customer routed to to this server is blocked and leaves the system. It should be clear that our proof
of Theorem~\ref{th-main} works for this more general system; additional ``losses'' of arriving customers
do not change the stochastic upper bound on a steady-state server workload; and the rest of the proof remains
essentially unchanged, except a more cumbersome state space description.

\bibliographystyle{acmtrans-ims}
\bibliography{biblio-stolyar}

\end{document}